\documentclass{article}
\date{}
\pdfoutput=1
\usepackage{hyperref}
\usepackage{indentfirst}
\usepackage{amsmath}
\usepackage{amssymb}
\usepackage{cases}
\title{Projective spherically symmetric Finsler metrics with constant flag curvature in $R^n$ }
\author{Linfeng Zhou}
\begin{document}
\maketitle
\begin{center}
\textit{Department of Mathematics}\\
\textit{East China Normal University, Shanghai, 200241, China} \\
\textit{E-mail: lfzhou@math.ecnu.edu.cn}
\end{center}

\begin{abstract}
We investigate projective spherically symmetric Finsler metrics with
constant flag curvature in $R^n$ and give the complete
classification theorems. Furthermore, a new class of Finsler metrics
with two parameters on n-dimensional disk are found to have constant
negative flag curvature.
\\

\noindent\textbf{2000 Mathematics Subject Classification:}
53B40, 53C60, 58B20.\\
\textbf{Keywords and Phrases: spherical symmetry, projective,
constant flag curvature, Bryant metrics.}
\end{abstract}

\newtheorem{Th}{Theorem}[section]
\newtheorem{Prop}[Th]{Propositon}
\newtheorem{Col}[Th]{Corollary}
\newtheorem{Lem}[Th]{Lemma}
\newtheorem{Ex}[Th]{Example}
\newtheorem{Con}[Th]{Conjecture}

\section{Introduction}
Suppose $\Omega\subseteq R^n$ is a convex domain and Finsler metric
$F$ is defined on $\Omega$. $(\Omega, F)$ is called spherically
symmetric if orthogonal matrix $O(n)$ is isometric map of $(\Omega,
F)$. It means that $(\Omega, F)$ is invariant under any rotations in
$R^n$. In \cite{Zhou}, using the Killing fields equation, the author
proves that Finsler metric $(\Omega,F)$ is spherically symmetric if
and only if $F$ can be written as $F=\phi(|x|,|y|,\langle
x,y\rangle)$.

Since this kind of metrics have a nice symmetry, many complex
computation can become much easier. Let $r=|x|$, $u=|y|$ and
$v=\langle x,y\rangle$, it is similar with $(\alpha,\beta)$ metric
to compute the fundamental tensor $g_{ij}$ of Finsler metric
$F=\phi(r,u,v) $\cite{Zhou}\cite{CS}:
\[g_{ij}=\frac{\phi\phi_u}{u}\delta_{ij}+(\phi_v^2+\phi\phi_{vv})x^ix^j+(\frac{\phi_u^2+\phi\phi_{uu}}{u^2}-\frac{\phi\phi_u}{u^3})y^iy^j+(\frac{\phi_u\phi_v+\phi\phi_{uv}}{u})(x^iy^j+x^jy^i).\]
Thus \label{eq18}
\begin{equation}det(g_{ij})=(\frac{\phi}{u})^{n+1}\phi_u^{n-2}[\phi_u+(r^2u^2-v^2)\frac{\phi_{vv}}{u}].\end{equation}

As we know, the Hilbert's Fourth problem relates to classify the
projective Finsler metric in $R^n$ which is still one of the
motivation in Finsler geometry \cite{Sh1}. A Finsler metric in $R^n$
is projective if it's geodesics are all straight lines. According to
Rapcs\'{a}k's lemma \cite{Ra}, a spherically symmetric Finsler
metric $F=\phi(r,u,v)$ is projective if and only if $\phi$ satisfies
that \cite{Zhou}
\begin{equation}\label{eq17}\phi_{rv}\frac{v}{r}+\phi_{vv}u^2=\frac{\phi_r}{r}.\end{equation}
Solving the equation can give a total classification theorem:
\begin{Th} \cite{Zhou}
Suppose $F$ is a spherically symmetric Finsler
metric on a convex domain $\Omega\in R^n$, $F$ is projective if and
only if there exist smooth functions $f(t)>0$ and $g(r)$ s.t.
\[\phi(r,u,v)=\int f(\frac{v^2}{u^2}-r^2)du+g(r)v\]
where $F(x,y)=\phi(|x|,|y|,\langle x,y\rangle)$.
\end{Th}

Furthermore, many classical Finsler metrics are projective spherically symmetric. We now list some of them here:
\begin{Ex}[Berwald metric \cite{Be}]
Let $B^n \subset R^n$ be a standard unit ball.
 An $(\alpha,\beta)$ metric $F$ is defined on $B^n$:
 \[F(x,y):=\frac{(\sqrt{|y|^2-(|x|^2|y|^2-\langle x,y\rangle^2)}+\langle x,y\rangle)^2}{(1-|x|^2)^2\sqrt{|y|^2-(|x|^2|y|^2-\langle x,y\rangle^2)}}.\]
 $F$ is projective and has constant flag curvature $K=0$.
 \end{Ex}

 \begin{Ex} [projective spherical model]
 Let $S^n \subset R^{n+1}$ be a standard unit sphere. The standard
 inner product $\langle,\rangle$ in $R^{n+1}$ induced a Riemannian
 metric on $S^n$: for $x\in S^n$, let
 \[\alpha:=|y|,\quad y\in T_xS^n\subset R^{n+1}.\]
Let $S^n_{+}$ denote the upper hemisphere and let
 $\psi_{+}:R^n\rightarrow S_{+}^n$ be the projection map defined by
 \[\psi_{+}(x):=(\frac{x}{\sqrt{1+|x|^2}},\frac{1}{1+|x|^2}).\]
The pull-back metric on $R^n$ from $S^n_{+}$ by $\psi_{+}$ is given
by
\[\alpha(x,y):=\frac{\sqrt{|y|^2+(|x|^2|y|^2-\langle x,y \rangle^2)}}{{1+|x|^2}},\quad y\in T_xR^n\]
$(R^n,\alpha(x,y))$ is projectively flat and has constant flag
curvature $K=1$.
\end{Ex}

\begin{Ex}[Bryant metric \cite{Br}]\label{ex1}
Denote
\[\begin{split}
A:=&(\cos(2\alpha)|y|^2+(|x|^2|y|^2-\langle x,y\rangle^2))^2+(\sin(2\alpha)|y|^2)^2,\\
B:=&\cos(2\alpha)|y|^2+(|x|^2|y|^2-\langle x,y\rangle^2),\\
C:=&\sin(2\alpha)\langle x,y\rangle,\\
D:=&|x|^4+2\cos(2\alpha)|x|^2+1.
\end{split}\]
For an angle $\alpha$ with $0\leq\alpha<\frac{\pi}{2}$, Bryant
metric $F$ is defined by
\[F:=\sqrt{\frac{\sqrt{A}+B}{2D}+(\frac{C}{D})^2}+\frac{C}{D}\]
on the whole region $R^n$. As we know it is projective and has
constant flag curvature $K=1$.
 \end{Ex}

 \begin{Ex}[Klein model] Let $B^n \subset R^n$ be a standard unit
 ball and let
 \[\alpha(x,y):=\frac{\sqrt{|y|^2-(|x|^2|y|^2-\langle x,y\rangle^2)}}{1-|x|^2},\quad y\in T_xB^n.\]
 $\alpha(x,y)$ is a Riemannnian metric on $B^n$. It is projective and
 has constant flag curvature $K=-1$.
 \end{Ex}

 \begin{Ex}[Funk metric \cite{BCS}]
 A Randers metric $F$ is defined on the standard unite ball $B^n$:
 \[F(x,y):=\frac{\sqrt{|y|^2-(|x|^2|y|^2-\langle x,y\rangle^2)}+\langle
 x,y\rangle}{1-|x|^2}.\] It is also projective and has constant flag
 curvature $K=-\frac{1}{4}$.
 \end{Ex}
Clearly, all of these examples have constant flag curvature. Hence a
natural question arises: are there any other projective spherically
symmetric Finsler metrics which also have constant flag curvature
and how to classify them?

In this paper, focusing on above question, we give an answer: by
solving the partial equations which characterize projective
symmetric Finsler metric in $R^n$ with constant flag curvature and
discussing the different cases, we obtain the complete
classification theorems \ref{th1}, \ref{th2} and \ref{th3}. In
theorem \ref{th3}, a new type of Finsler metrics with constant flag
curvature are discovered.

\section{Characterizing Equations}
In \cite{Zhou}, the author gives the following theorem which
contains two partial equations characterizing projective Finsler
metric $F=\phi(|x|,|y|,\langle x,y\rangle)$ with constant flag
curvature:
\begin{Th}\label{th4} Let a spherically symmetric Finsler metric $F=\phi(|x|,|y|,\langle
x,y\rangle)$ be projective on a convex domain $\Omega \subseteq
R^n$. Then $F$ has constant flag curvature $K=\lambda$ if and only
if $\phi(r,u,v)$ satisfies
\begin{subnumcases}
{}4\lambda r\phi^4\phi_u+r\phi_uQ^2-4ru\phi\phi_vQ+4u\phi^2\phi_r=0\label{eq:a1}\\
4\lambda r\phi^4\phi_v+r\phi_vQ^2+2\phi^2
          Q_r-4\phi\phi_rQ=0\label{eq:a2}
\end{subnumcases}
where $Q:=\frac{v}{r}\phi_r+u^2\phi_v$.
\end{Th}

In order to solve the equations, we need change variables to
simplify them. Let

\begin{equation*}
\left\{ \begin{aligned}
         r&=\sqrt{z_1^2+z_2^2}\\
         v&=uz_2\\
         u&=u.
                          \end{aligned} \right.
                          \end{equation*}
Under the new parameters, assume
$F=\phi(r,u,v)=\overline{\phi}(z_1,z_2,u)$. Hence we compute
\begin{equation} \label{eq2}
\left\{ \begin{aligned}
         \phi_u&=\overline{\phi}_{z_1}\frac{z_2^2}{uz_1}-\overline{\phi}_{z_2}\frac{z_2}{u}+\overline{\phi}_u\\
         \phi_v&=-\overline{\phi}_{z_1}\frac{z_2}{uz_1}+\overline{\phi}_{z_2}\frac{1}{u}\\
         \phi_r&=\overline{\phi}_{z_1}\frac{r}{z_1}
                          \end{aligned} \right.
                          \end{equation}
and
\begin{equation}\label{eq3}Q=\frac{v}{r}\phi_r+u^2\phi_v=u\overline{\phi}_{z_2}.\end{equation}
Plugging (\ref{eq2}) and (\ref{eq3}) into the equation
(\ref{eq:a1}), we obtain
\begin{align*}\frac{\sqrt{z_1^2+z_2^2}}{uz_1}(&4\lambda\overline{\phi}^4\overline{\phi}_{z_1}z_2^2-4\lambda\overline{\phi}^4\overline{\phi}_{z_2}z_1z_2+4\lambda\overline{\phi}^5z_1+u^2\overline{\phi}_{z_2}^2\overline{\phi}_{z_1}z_2^2-u^2\overline{\phi}_{z_2}^3z_1z_2\\
&-3u^2\overline{\phi}_{z_2}^2\overline{\phi}z_1+4u^2\overline{\phi}_{z_2}\overline{\phi}_{z_1}\overline{\phi}z_2+4u^2\overline{\phi}^2\overline{\phi}_{z_1})=0.\end{align*}
So
\[\overline{\phi}_{z_1}=\frac{z_1(u^2\overline{\phi}_{z_2}^3z_2+4\lambda\overline{\phi}^4\overline{\phi}_{z_2}z_2-4\lambda\overline{\phi}^5+3u^2\overline{\phi}_{z_2}^2\overline{\phi})}{4\lambda\overline{\phi}^4z_2^2
+u^2\overline{\phi}_{z_2}^2z_2^2+4u^2\overline{\phi}_{z_2}\overline{\phi}z_2+4u^2\overline{\phi}^2}.\]
Observing the homogeneity of Finsler metric
$F=\phi(r,u,v)=\overline{\phi}(z_1,z_2,u)$, we have
\[F=\overline{\phi}(z_1,z_2,u)=u\widetilde{\phi}(z_1,z_2).\]
Thus the above equation can be written as
\begin{equation}\label{eq4}
\widetilde{\phi}_{z_1}=\frac{z_1(\widetilde{\phi}_{z_2}^3z_2+4\lambda\widetilde{\phi}^4\overline{\phi}_{z_2}z_2-4\lambda\widetilde{\phi}^5+3\widetilde{\phi}_{z_2}^2\widetilde{\phi})}{4\lambda\widetilde{\phi}^4z_2^2
+\widetilde{\phi}_{z_2}^2z_2^2+4\widetilde{\phi}_{z_2}\widetilde{\phi}z_2+4\widetilde{\phi}^2}.
\end{equation}
Similarly, substituting (\ref{eq2}), (\ref{eq3}) and
$Q_r=\frac{u\sqrt{z_1^2+z_2^2}}{z_1}\overline{\phi}_{z_1z_2}$ to the
equation (\ref{eq:a2}) yields the following equation
\begin{align*}\sqrt{z_1^2+z_2^2}\overline{\phi}(&16\lambda^2\overline{\phi}^8z_1z_2-8\lambda\overline{\phi}^4z_1u^2\overline{\phi}_{z_2}^2z_2+32\lambda\overline{\phi}^5z_1u^2\overline{\phi}_{z_2}-3u^4\overline{\phi}_{z_2}^4z_1z_2-8u^4\overline{\phi}_{z_2}^3z_1\overline{\phi}\\
&+8\overline{\phi}^5u^2\overline{\phi}_{z_1z_2}\lambda
z_2^2+2\overline{\phi}
u^4\overline{\phi}_{z_1z_2}\overline{\phi}_{z_2}^2z_2^2+8\overline{\phi}^2u^4\overline{\phi}_{z_1z_2}\overline{\phi}_{z_2}z_2+8\overline{\phi}^3u^4\overline{\phi}_{z_1z_2})\\
&\div[(4\lambda\overline{\phi}^4z_2^2+u^2\overline{\phi}_{z_2}^2z_2^2+4u^2\overline{\phi}_{z_2}\overline{\phi}z_2+4u^2\overline{\phi}^2)uz_1]=0.\end{align*}
Therefore
\[\overline{\phi}_{z_1z_2}=\frac{z_1(-16\lambda^2\overline{\phi}^8z_2+8\lambda u^2\overline{\phi}^4\overline{\phi}_{z_2}^2z_2-32\lambda u^2\overline{\phi}^5\overline{\phi}_{z_2}+3u^4\overline{\phi}_{z_2}^4z_2+8u^4\overline{\phi}_{z_2}^3\overline{\phi})}{2u^2\overline{\phi}(4\lambda\overline{\phi}^4z_2^2+u^2\overline{\phi}_{z_2}^2z_2^2+4u^2\overline{\phi}_{z_2}\overline{\phi}z_2+4u^2\overline{\phi}^2)}.\]
Since $F=\overline{\phi}=u\widetilde{\phi}(z_1,z_2)$, the above
equation implies
\begin{equation} \label{eq5}
\widetilde{\phi}_{z_1z_2}=\frac{z_1(-16\lambda^2\widetilde{\phi}^8z_2+8\lambda
\widetilde{\phi}^4\widetilde{\phi}_{z_2}^2z_2-32\lambda
\widetilde{\phi}^5\widetilde{\phi}_{z_2}+
3\widetilde{\phi}_{z_2}^4z_2+8\widetilde{\phi}_{z_2}^3\widetilde{\phi})}
{2\widetilde{\phi}(4\lambda\widetilde{\phi}^4z_2^2+\widetilde{\phi}_{z_2}^2z_2^2
+4\widetilde{\phi}_{z_2}\widetilde{\phi}z_2+4\widetilde{\phi}^2)}.
\end{equation}
Combining (\ref{eq4})-(\ref{eq5}) and
$(\widetilde{\phi}_{z_1})_{z_2}=\widetilde{\phi}_{z_1z_2}$, we find
\begin{align*}
(3\widetilde{\phi}_{z_2}^2-&2\widetilde{\phi}\widetilde{\phi}_{z_2z_2}-4\lambda\widetilde{\phi}^4)(8\widetilde{\phi}\widetilde{\phi}_{z_2}^3z_2^2+24\widetilde{\phi}^2\widetilde{\phi}_{z_2}^2z_2+\widetilde{\phi}_{z_2}^4z_2^3+24\widetilde{\phi}^3\widetilde{\phi}_{z_2}\\
&+16\lambda^2z_2^3\widetilde{\phi}^8+32\lambda
z_2^2\widetilde{\phi}^5\widetilde{\phi}_{z_2}+32\lambda
z_2\widetilde{\phi}^6+8\lambda
z_2^3\widetilde{\phi}^4\widetilde{\phi}_{z_2}^2)=0.\end{align*} It
means that either
\begin{equation} \label{eq6}
3\widetilde{\phi}_{z_2}^2-2\widetilde{\phi}\widetilde{\phi}_{z_2z_2}-4\lambda\widetilde{\phi}^4=0
\end{equation}
or
\begin{equation}\label{eq7}
\begin{aligned}
&\widetilde{\phi}_{z_2}(2\widetilde{\phi}+z_2\widetilde{\phi}_{z_2})(z_2^2\widetilde{\phi}_{z_2}^2+6z_2\widetilde{\phi}\widetilde{\phi}_{z_2}+12\widetilde{\phi}^2)\\
&+8\lambda z_2\widetilde{\phi}^4(2\lambda
z_2^2\widetilde{\phi}^4+4\widetilde{\phi}^2+4z_2\widetilde{\phi}\widetilde{\phi}_{z_2}+z_2^2\widetilde{\phi}_{z_2}^2)=0
\end{aligned}
\end{equation}
holds. For the equation (\ref{eq6}), we have following lemma.
\begin{Lem}\label{lem1}
The solutions of ODE
\[2y\frac{d^2y}{dx^2}-3(\frac{dy}{dx})^2+4\lambda y^4=0\]
are
\[y=0\quad\textit{and}\quad y=\frac{1}{c_1x^2+c_2x+c_3} \]
where $c_1$, $c_2$, $c_3$ are constant and satisfy
\[c_2^2-4c_1c_3+4\lambda=0.\]
\end{Lem}
\emph{Proof.} Obviously $y=0$ solves the ODE. If $y\neq0$, use
$w=\frac{1}{y}$ to replace $y$ and compute
\[\frac{dy}{dx}=-\frac{1}{w^2}\frac{dw}{dx}\quad \textit{and}\quad \frac{d^2y}{dx^2}=-\frac{1}{w^2}\frac{d^2w}{dx^2}+\frac{2}{w^3}(\frac{dw}{dx})^2.\]
Hence the ODE becomes
\begin{equation}\label{eq8}(\frac{dw}{dx})^2-2w\frac{d^2w}{dx^2}+4\lambda=0.\end{equation} Derivative
of above equation yields
\[2w\frac{d^3w}{dx^3}=0.\]
Since $w\neq0$, so
\[\frac{d^3w}{dx^3}=0.\]
Thus
\[w=c_1x^2+c_2x+c_3.\]
Plugging it into (\ref{eq8}), we have
\[(2c_1x+c_2)^2-2(c_1x^2+c_2x+c_3)2c_1+4\lambda=0.\]
This implies that $c_1$, $c_2$ and $c_3$ satisfy
\[c_2^2-4c_1c_3+4\lambda=0,\]
as asserted.

Conversely, it is easy to verify that $y=\frac{1}{c_1x^2+c_2x+c_3}$
solves the ODE if $c_1$, $c_2$, $c_3$ satisfy
$c_2^2-4c_1c_3+4\lambda=0$.\qquad Q.E.D.\\

For simplicity, we may assume that constant flag curvature
$K=0,1,-1$. For other case, by a scaling transformation, it can be
changed into these three cases. Now we are ready to prove the
classification theorems of projective spherically symmetric Finsler
metrics with constant flag curvature in $R^n$.

\section{Vanishing flag curvature}

Let us firstly consider the case of flag curvature $K=0$:
\begin{Th} \label{th1}
On a convex domain $\Omega \subseteq R^n$, a projective spherically
symmetric Finsler metric $F=\phi(|x|,|y|,\langle x,y\rangle)$ has
constant flag curvature $K=0$ if and only if either
\begin{enumerate}
\item[(1)] $\Omega=R^n$ with Euclidean metric $F=|y|$; or
\item[(2)] $\Omega=B^n(\sqrt{c}):=\{x:|x|^2<c\}$ with Berwald metric
\[F=\frac{|y|}{\sqrt{c-z_1^2}(z_2\pm\sqrt{c-z_1^2})^2}\]
where $z_1:=\sqrt{|x|^2-\frac{\langle x,y\rangle^2}{|y|^2}},
z_2:=\frac{\langle x,y\rangle}{|y|}$.
\end{enumerate}
\end{Th}
\emph{Proof.} Sufficiency is obvious. Now we only need to prove the
necessity. According to the analysis in last section, the projective
Finsler metric $F=u\widetilde{\phi}(z_1,z_2)$ which has constant
flag curvature must satisfy either (\ref{eq6}) or (\ref{eq7}). When
$K=\lambda=0$, (\ref{eq6}) and (\ref{eq7}) can be written as
\begin{align}\label{eq6:1}
3\widetilde{\phi}_{z_2}^2-2\widetilde{\phi}\widetilde{\phi}_{z_2z_2}=0\quad
\textit{and} \tag{\ref{eq6}$'$}
\end{align}
\begin{align}\label{eq7:1}
\widetilde{\phi}_{z_2}(2\widetilde{\phi}+z_2\widetilde{\phi}_{z_2})(z_2^2\widetilde{\phi}_{z_2}^2+6z_2\widetilde{\phi}\widetilde{\phi}_{z_2}+12\widetilde{\phi}^2)=0.
\tag{\ref{eq7}$'$}
\end{align}
If (\ref{eq7:1}) holds, it is equivalent to have either
\[\widetilde{\phi}_{z_2}=0\quad \textit{or}\]
\[2\widetilde{\phi}+z_2\widetilde{\phi}_{z_2}=0\quad \textit{or}\]
\[z_2^2\widetilde{\phi}_{z_2}^2+6z_2\widetilde{\phi}\widetilde{\phi}_{z_2}+12\widetilde{\phi}^2=0.\]
\begin{enumerate}
\item[(a)]When $\widetilde{\phi}_{z_2}=0$, from (\ref{eq4}) we find
\[\widetilde{\phi}_{z_1}=0.\]
It means that Finsler metric is Euclidean metric and $F=cu$,
$\Omega=R^n$. By a scaling, we may let $F=u$.
\item [(b)]When
$2\widetilde{\phi}+z_2\widetilde{\phi}_{z_2}=0$, it can be solved
that $\widetilde{\phi}=\frac{c(z_1)}{z_2^2}$. Hence
$F=u\widetilde{\phi}=u\frac{c(z_1)}{z_2^2}$. Recalling
$z_2=\frac{v}{u}=\frac{\langle x,y\rangle}{|y|}$, we can observe
that $F$ is singular on some directions. Consequently, this kind of
situation is impossible.
\item [(c)]When
$z_2^2\widetilde{\phi}_{z_2}^2+6z_2\widetilde{\phi}\widetilde{\phi}_{z_2}+12\widetilde{\phi}^2=0$,
solving $\widetilde{\phi}_{z_2}$ from it, we obtain
\[\widetilde{\phi_{z_2}}=\frac{-6z_2\pm\sqrt{-12z_2^2\widetilde{\phi}^2}}{2z_2^2}.\]
Clearly, there are no real solutions.
\end{enumerate}
Now we mainly focus on studying (\ref{eq6:1}). By lemma \ref{lem1},
we know
\[\widetilde{\phi}(z_1,z_2)=\frac{1}{\widetilde{c}_1(z_1)z_2^2+\widetilde{c}_2(z_1)z_2+\widetilde{c}_3(z_1)}\]
where $\widetilde{c}_2^2-4\widetilde{c}_1\widetilde{c}_3=0$.
\begin{enumerate}
\item[(a)] If $\widetilde{c}_1=0$, then $\widetilde{c}_2=0$ and $\widetilde{\phi}=\frac{1}{\widetilde{c}_3(z_1)}$.
Therefore we have $\widetilde{\phi}_{z_2}=0$. Similar with above
analysis, we conclude that $F$ is Euclidean metric and $\Omega=R^n$.
\item[(b)] If $\widetilde{c}_1\neq0$, then $\widetilde{c}_3=\frac{\widetilde{c}_2^2}{4\widetilde{c}_1}$. So
\begin{align*}\widetilde{\phi}(z_1,z_2)&=\frac{1}{\widetilde{c}_1z_2^2+\widetilde{c}_2z_2+\frac{\widetilde{c}_2^2}{4\widetilde{c}_1}}\\
&=\frac{1}{\widetilde{c}_1(z_2+\frac{\widetilde{c}_2}{2\widetilde{c}_1})^2}=\frac{1}{c_1(z_1)(z_2+c_2(z_1))^2}.\end{align*}
Here we denote $c_1=\widetilde{c}_1$ and
$c_2=\frac{\widetilde{c}_2}{2\widetilde{c}_1}$. Plugging
$\widetilde{\phi}$ into the equation (\ref{eq4}) and computing
yields
\[(z_1c_1+c_2^2c_1')z_2+3z_1c_1c_2+c_2^3c_1'+2c_1c_2^2c_2'=0.\]
Since it holds for all $z_2\in R$, there must have
\begin{equation*}
\left\{ \begin{aligned} &z_1c_1+c_2^2c_1'=0\\
&3z_1c_1c_2+c_2^3c_1'+2c_1c_2^2c_2'=0.
                          \end{aligned} \right.
                          \end{equation*}
From the first equation, we obtain $c_2\neq0$ and
$\frac{c_1'c_2^2}{c_1}=-z_1$. By the second equation, we find that
$3z_1+\frac{c_1'c_2^2}{c_1}+2c_2c_2'=0$. Combining these two
identities yields
\[c_2c_2'=-z_1.\]
So we can solve that
\[c_2=\pm\sqrt{c-z_1^2}\quad \textit{and}\quad c_1=d\sqrt{c-z_1^2}\]
where $c$, $d$ are positive constants. Thereby,
\[F=\overline{\phi}(z_1,z_2,u)=u\widetilde{\phi}(z_1,z_2)=\frac{u}{d\sqrt{c-z_1^2}(z_2\pm\sqrt{c-z_1^2})^2}.\]
By a scaling, we can let
\[F=\frac{u}{\sqrt{c-z_1^2}(z_2\pm\sqrt{c-z_1^2})^2}\]
as asserted. Obviously, $F$ can be defined on
$\Omega=B^n(\sqrt{c})$.\qquad Q.E.D.
\end{enumerate}
\noindent\emph{Remark.} When $c=1$ in (2) of theorem \ref{th1},
Finsler metric $F$ becomes
\begin{align*}F&=\frac{|y|}{\sqrt{1-z_1^2}(z_2\pm\sqrt{1-z_1^2})^2}=\frac{|y|(\sqrt{1-z_1^2}\mp z_2)^2}{\sqrt{1-z_1^2}(1-z_1^2-z_2^2)^2}\\
&=\frac{(\sqrt{|y|^2-(|x|^2|y|^2-\langle x,y\rangle^2)}\mp\langle
x,y\rangle)^2}{(1-|x|^2)^2\sqrt{|y|^2-(|x|^2|y|^2-\langle
x,y\rangle^2)}}.\end{align*} Thus it is Berwald metric.

\section{Constant positive flag curvature}
Now we discuss the case of flag curvature $K=1$:
\begin{Th} \label{th2}
On a convex domain $\Omega \subseteq R^n$, a projective spherically
symmetric Finsler metric $F=\phi(|x|,|y|,\langle x,y\rangle)$ has
constant flag curvature $K=1$ if and only if either
\begin{enumerate}
\item[(1)] $\Omega=R^n$ with projective spherical metric
\[F=\frac{\sqrt{|y|^2+(|x|^2|y|^2-\langle x,y \rangle^2)}}{{1+|x|^2}}; \textit{or}\]
\item[(2)] $\Omega=R^n$ with Bryant type metric
\[F=\frac{|y|c_1(z_1)}{c_1(z_1)^2+\big(z_2+c_2(z_1)\big)^2}\] where
\[z_1:=\sqrt{|x|^2-\frac{\langle x,y\rangle^2}{|y|^2}}, z_2:=\frac{\langle x,y\rangle}{|y|},\]
\[c_1(z_1):=\frac{\sqrt{2}}{2}\sqrt{2d_2+z_1^2+\sqrt{(2d_2+z_1^2)^2+4d_1^2}},\]
\[c_2(z_1):=\pm\frac{\sqrt{2}}{2}\sqrt{-2d_2-z_1^2+\sqrt{(2d_2+z_1^2)^2+4d_1^2}},\]$d_1$
and $d_2$ are positive constants.
\end{enumerate}
\end{Th}
\emph{Proof.} Also, we only need to prove the necessity. Under the
flag curvature condition $K=\lambda=1$, the equations (\ref{eq6})
and (\ref{eq7}) come to
\begin{align}\label{eq6:2}
3\widetilde{\phi}_{z_2}^2-2\widetilde{\phi}\widetilde{\phi}_{z_2z_2}-4\widetilde{\phi}^4=0\quad
\textit{and} \tag{\ref{eq6}$''$}
\end{align}
\begin{align}\label{eq7:2}
&\widetilde{\phi}_{z_2}(2\widetilde{\phi}+z_2\widetilde{\phi}_{z_2})(z_2^2\widetilde{\phi}_{z_2}^2+6z_2\widetilde{\phi}\widetilde{\phi}_{z_2}+12\widetilde{\phi}^2)\nonumber\\
&+8z_2\widetilde{\phi}^4(2
z_2^2\widetilde{\phi}^4+4\widetilde{\phi}^2+4z_2\widetilde{\phi}\widetilde{\phi}_{z_2}+z_2^2\widetilde{\phi}_{z_2}^2)=0.
\tag{\ref{eq7}$''$}
\end{align}
If (\ref{eq7:2}) holds, solving it by Maple we know that
\[\widetilde{\phi}=\frac{f(z_1,z_2)}{z_2}.\]
Hence $F=u\widetilde{\phi}$ is also singular on some directions.

Now we consider the equation (\ref{eq6:2}). Due to lemma \ref{lem1},
we know
\[\widetilde{\phi}(z_1,z_2)=\frac{1}{\widetilde{c}_1(z_1)z_2^2+\widetilde{c}_2(z_1)z_2+\widetilde{c}_3(z_1)}\]
where $\widetilde{c}_2^2-4\widetilde{c}_1\widetilde{c}_3+4=0$.
Obviously, $\widetilde{c}_1\neq0$, then
$\frac{\widetilde{c}_3}{\widetilde{c}_1}=\frac{1}{\widetilde{c}_1^2}+\frac{\widetilde{c}_2^2}{4\widetilde{c}_1^2}$.
Substituting it to $\widetilde{\phi}$ and letting
$c_1=\frac{1}{\widetilde{c}_1}$,
$c_2=\frac{\widetilde{c}_2}{2\widetilde{c}_1}$, we have
\[\widetilde{\phi}=\frac{c_1(z_1)}{c_1(z_1)^2+(z_2+c_2(z_1))^2}.\]
Plugging the formula of $\widetilde{\phi}$ into the equation
(\ref{eq4}) yields
\begin{align*}(z_1c_1-c_1'c_2^2-&c_1'c_1^2)z_2^2+(-2c_1'c_1^2c_2+2c_1c_2^2c_2'-2c_1'c_2^3+2c_1^3c_2'+4z_1c_1c_2)z_2\\
&-c_1'c_2^4+2c_1c_2^3c_2'+2c_1^3c_2c_2'+3z_1c_1c_2^2+c_1^4c_1'-z_1c_1^3=0.\end{align*}
It holds for all $z_2\in R$ if and only if
\begin{equation}\label{eq11}
\left\{ \begin{aligned} &z_1c_1-c_1'c_2^2-c_1'c_1^2=0\\
&-2c_1'c_1^2c_2+2c_1c_2^2c_2'-2c_1'c_2^3+2c_1^3c_2'+4z_1c_1c_2=0\\
&-c_1'c_2^4+2c_1c_2^3c_2'+2c_1^3c_2c_2'+3z_1c_1c_2^2+c_1^4c_1'-z_1c_1^3=0.
                          \end{aligned} \right.
                          \end{equation}
Solving above equations and getting rid of the improper solutions
will get
\begin{equation} \label{eq9}
\left\{ \begin{aligned} &c_1(z_1)=\sqrt{z_1^2+c}\\
&c_2(z_1)=0
                          \end{aligned} \right.
                          \end{equation}
and
\begin{equation}\label{eq10}
\left\{ \begin{aligned} &c_1(z_1)=\frac{\sqrt{2}}{2}\sqrt{2d_2+z_1^2+\sqrt{(2d_2+z_1^2)^2+4d_1^2}}\\
&c_2(z_1)=\pm\frac{\sqrt{2}}{2}\sqrt{-2d_2-z_1^2+\sqrt{(2d_2+z_1^2)^2+4d_1^2}}.
                          \end{aligned} \right.
                          \end{equation}
See the following lemma for more details.

In the case of (\ref{eq9}), Finsler metric
$F=\frac{u\sqrt{c+z_1^2}}{c+z_1^2+z_2^2}=\frac{\sqrt{c|y|^2+|x|^2|y|^2-\langle
x,y\rangle^2}}{c+|x|^2}$ is projective spherical metric which can be
defined on $R^n$.

In the case of (\ref{eq10}), Finsler metric
$F=\frac{uc_1}{c_1^2+(z_2+c_2)^2}$ is Bryant type metric which can
also be defined on $R^n$.\qquad Q.E.D.\\

\noindent\emph{Remark.} It can be verified by Maple: when
$d_1=\frac{1}{4}\sin^2(2\alpha)$, $d_2=\frac{1}{2}\cos(2\alpha)$,
Finsler metric in (2) of theorem \ref{th2} turns to Bryant metric in
example \ref{ex1}.

\begin{Lem} The solutions of equation system (\ref{eq11}) are

\begin{equation*}
\left\{ \begin{aligned} &c_1=0\\
&c_2=c_2,
                          \end{aligned} \right.
                          \quad\left\{
                          \begin{aligned}&c_1=\pm\sqrt{z_1^2+c}\\
                          &c_2=0
                          \end{aligned}\right.
                          \end{equation*}
and
\begin{equation*}
\left\{ \begin{aligned} &c_1=\pm\frac{\sqrt{2}}{2}\sqrt{2d_2+z_1^2\pm\sqrt{(z_1^2+2d_2)^2+4d_1^2}}\\
&c_2=\pm\frac{\sqrt{2}}{2}\sqrt{-2d_2-z_1^2\pm\sqrt{(z_1^2+2d_2)^2+4d_1^2}}.
                          \end{aligned} \right.
                          \end{equation*}
\end{Lem}
\emph{Proof.} The equation system (\ref{eq11}) can be rewritten as
\begin{equation}\label{eq12}
\left\{ \begin{aligned} &c_1'(c_1^2+c_2^2)=z_1c_1\\
&c_1'c_2(c_1^2+c_2^2)-c_1c_2'(c_1^2+c_2^2)-2z_1c_1c_2=0\\
&c_1'(c_1^4-c_2^4)+2c_2'c_1c_2(c_1^2+c_2^2)+3z_1c_1c_2^2-z_1c_1^3=0.
                          \end{aligned} \right.
                          \end{equation}
Obviously, $c_1=0$, $c_2=c_2$ is one of the solutions. Now assuming
$c_1\neq0$, from the first equation of (\ref{eq12}) we have
\begin{equation}\label{eq13}c_1'=\frac{z_1c_1}{c_1^2+c_2^2}.\end{equation}
Plugging it into the second equation of (\ref{eq12}), we will get
\begin{equation}\label{eq14}c_2'=-\frac{z_1c_2}{c_1^2+c_2^2}.\end{equation}
Substituting (\ref{eq13}) and (\ref{eq14}) to the third equation of
(\ref{eq12}) will find it automatically holds. Thus the original
equation system is equivalent to the equation (\ref{eq13}) and
(\ref{eq14}). Together with these two equations we know
\[(c_1c_2)'=c_1'c_2+c_2'c_1=0.\]
That means there exists a constant $d_1$ s.t.
\[c_2=\frac{d_1}{c_1}.\]
If $d_1=0$, $c_2=0$ and (\ref{eq13}) becomes
\[c_1'c_1=z_1.\]
Thus $c_1=\pm\sqrt{z_1^2+c}$ where $c$ is a constant.
 Otherwise, (\ref{eq13}) becomes
\[c_1'(c_1+\frac{d_1^2}{c_1^2})=z_1.\]
Solve it by an integration to find
\[c_1=\pm\frac{\sqrt{2}}{2}\sqrt{2d_2+z_1^2\pm\sqrt{(z_1^2+2d_2)^2+4d_1^2}}\]
where $d_2$ is a constant. So
\[c_2=\pm\frac{\sqrt{2}}{2}\sqrt{-2d_2-z_1^2\pm\sqrt{(z_1^2+2d_2)^2+4d_1^2}}.\]
Q.E.D.\\
\section{Constant negative flag curvature}
Finally, we study the case of flag curvature $K=-1$.
\begin{Th} \label{th3}
On a convex domain $\Omega \subseteq R^n$, a projective spherically
symmetric Finsler metric $F=\phi(|x|,|y|,\langle x,y\rangle)$ has
constant flag curvature $K=-1$ if and only if either
\begin{enumerate}
\item[(1)] $\Omega=B^n(\sqrt{c})$ with Klein metric
\[F=\frac{\sqrt{c|y|^2-(|x|^2|y|^2-\langle x,y\rangle^2)}}{c-|x|^2}; \textit{or}\]
\item[(2)] $\Omega=B^n(\sqrt{c})$ with Randers metric
\[F=\frac{\sqrt{c|y|^2-(|x|^2|y|^2-\langle x,y\rangle^2)}+\langle x,y\rangle}{2(c-|x|^2)}; \textit{or}\]
\item[(3)] $\Omega=B^n(\sqrt{2(d_2-d_1)})$ with Finsler metric
\[F=\frac{|y|c_1(z_1)}{c_1(z_1)^2-\big(z_2+c_2(z_1)\big)^2}\] where
\[z_1:=\sqrt{|x|^2-\frac{\langle x,y\rangle^2}{|y|^2}}, z_2:=\frac{\langle x,y\rangle}{|y|},\]
\[c_1(z_1):=\frac{\sqrt{2}}{2}\sqrt{2d_2-z_1^2+\sqrt{(2d_2-z_1^2)^2-4d_1^2}},\]
\[c_2(z_1):=\pm\frac{\sqrt{2}}{2}\sqrt{2d_2-z_1^2-\sqrt{(2d_2-z_1^2)^2-4d_1^2}},\]$d_1$
and $d_2$ are positive constant.
\end{enumerate}
\end{Th}
\emph{Proof.} Indicated by (\ref{eq18}), (\ref{eq17}) and
characterizing equations in theorem \ref{th4}, one can prove the
sufficiency through a direct calculation. Now we prove the
necessity. When the flag curvature $K=\lambda=-1$, the equations
(\ref{eq6}) and (\ref{eq7}) are
\begin{align}\label{eq6:3}
3\widetilde{\phi}_{z_2}^2-2\widetilde{\phi}\widetilde{\phi}_{z_2z_2}+4\widetilde{\phi}^4=0\quad
\textit{and} \tag{\ref{eq6}$'''$}
\end{align}
\begin{align}\label{eq7:3}
&\widetilde{\phi}_{z_2}(2\widetilde{\phi}+z_2\widetilde{\phi}_{z_2})(z_2^2\widetilde{\phi}_{z_2}^2+6z_2\widetilde{\phi}\widetilde{\phi}_{z_2}+12\widetilde{\phi}^2)\nonumber\\
&-8z_2\widetilde{\phi}^4(2
z_2^2\widetilde{\phi}^4+4\widetilde{\phi}^2+4z_2\widetilde{\phi}\widetilde{\phi}_{z_2}+z_2^2\widetilde{\phi}_{z_2}^2)=0.
\tag{\ref{eq7}$'''$}
\end{align}
Analogously, Finsler metric can't satisfy (\ref{eq7:3}). Thus we
only need consider the equation (\ref{eq6:3}). Because of lemma
\ref{lem1} and some computation, we have
\[\widetilde{\phi}=\frac{1}{2(z_2+c_1(z_1))}\quad \textit{or}\quad \widetilde{\phi}=\frac{c_1(z_1)}{c_1(z_1)^2-(z_2+c_2(z_1))^2} .\]

If $\widetilde{\phi}=\frac{1}{2(z_2+c_1(z_1))}$, plugging it into
the equation (\ref{eq4}) obtains
\[c_1'c_1+z_1=0.\]
So $c_1=\pm\sqrt{c-z_1^2}$ and $F$ becomes
\begin{align*}F=\phi&=\frac{u}{2(z_2\pm\sqrt{c-z_1^2})}=\frac{u(\sqrt{c-z_1^2}\mp z_2)}{2(c-z_1^2-z_2^2)}\\
&=\frac{\sqrt{c|y|^2-(|x|^2|y|^2-\langle x,y\rangle^2)}\mp\langle
x,y\rangle}{2(c-|x|^2)}.\end{align*} That is a Randers metric which
is defined on $B^n(\sqrt{c})$.

If $\widetilde{\phi}=\frac{c_1(z_1)}{c_1(z_1)^2-(z_2+c_2(z_1))^2}$,
from the equation (\ref{eq4}) we have
\begin{align*}(-z_1c_1+c_1'c_2^2-&c_1'c_1^2)z_2^2+(-2c_1'c_1^2c_2-2c_1c_2^2c_2'+2c_1'c_2^3+2c_1^3c_2'-4z_1c_1c_2)z_2\\
&+c_1'c_2^4-2c_1c_2^3c_2'+2c_1^3c_2c_2'-3z_1c_1c_2^2-c_1^4c_1'-z_1c_1^3=0.\end{align*}
It holds for all $z_2\in R$ if and only if
\begin{equation*}
\left\{ \begin{aligned} &-z_1c_1+c_1'c_2^2-c_1'c_1^2=0\\
&-2c_1'c_1^2c_2-2c_1c_2^2c_2'+2c_1'c_2^3+2c_1^3c_2'-4z_1c_1c_2=0\\
&c_1'c_2^4-2c_1c_2^3c_2'+2c_1^3c_2c_2'-3z_1c_1c_2^2-c_1^4c_1'-z_1c_1^3=0.
                          \end{aligned} \right.
                          \end{equation*}
Solving above equations and discarding the unsuitable solutions will
yield
\begin{equation} \label{eq15}
\left\{ \begin{aligned} &c_1(z_1)=\sqrt{c-z_1^2}\\
&c_2(z_1)=0
                          \end{aligned} \right.
                          \end{equation}
and
\begin{equation}\label{eq16}
\left\{ \begin{aligned} &c_1(z_1)=\frac{\sqrt{2}}{2}\sqrt{2d_2-z_1^2+\sqrt{(2d_2-z_1^2)^2-4d_1^2}}\\
&c_2(z_1)=\pm\frac{\sqrt{2}}{2}\sqrt{2d_2-z_1^2-\sqrt{(2d_2-z_1^2)^2-4d_1^2}}.
                          \end{aligned} \right.
                          \end{equation}

In the case of (\ref{eq15}), Finsler metric
$F=\frac{u\sqrt{c-z_1^2}}{c-z_1^2-z_2^2}=\frac{\sqrt{c|y|^2-(|x|^2|y|^2-\langle
x,y\rangle^2)}}{c-|x|^2}$ is Klein metric which is defined on
$B^n(\sqrt{c})$.

In the case of (\ref{eq16}), Finsler metric
$F=\frac{uc_1}{c_1^2-(z_2+c_2)^2}$ is a new Finsler metric which can
be defined on $B^n(\sqrt{2(d_2-d_1)})$.\qquad Q.E.D.\\

\noindent\emph{Remark.} To some extend, like Klein model and
projective spherical model, Finsler metric defined in (3) of theorem
\ref{th3} can be viewed as one of the pair of Bryant metric.

 \LaTeX
\end{document}